\documentstyle[12pt]{article}
\begin{document}

\title{Higher $l$-adic Abel-Jacobi mappings and filtrations
on Chow groups}
\author{Wayne Raskind\thanks{Research supported by the National Science
Foundation and the National Security Agency}\\Dept. of Mathematics\\University
of
Southern California \\
Los Angeles, CA 90089}
\date{}
\maketitle

\newcommand{\kb}{{\overline k}}
\newcommand{\bQ}{\mbox{\bf Q}}
\newcommand{\bZ}{\mbox{\bf Z}}
\newcommand{\bQl}{\mbox{\bf Q}_l}
\newcommand{\bZl}{\mbox{\bf Z}_l}
\newcommand{\sX}{{\cal X}}
\newcommand{\sS}{{\cal S}}
\newcommand{\sO}{{\cal O}}
\newcommand{\noi}{\noindent}
\newcommand{\plim}{\lim_{\stackrel{\leftarrow}{m}}}
\newcommand{\dlim}{\lim_{\stackrel{\rightarrow}{n}}}
\newcommand{\Xb}{{\overline X}}
\newcommand{\cF}{{\cal F}}
\newcommand{\cC}{{\cal C}}
\newcommand{\bF}{{\bf F}}
\newcommand{\cK}{{\cal K}}
\newcommand{\sU}{{\cal U}}
\newcommand{\sY}{{\cal Y}}

\section*{Introduction}

	Let $X$ be a smooth projective variety over a field
$k$ and denote by
$CH^n(X)$ the group of codimension $n$ cycles modulo
rational equivalence on $X$.
The structure of this group is largely unknown when $n\geq
2$, especially when
$k$ is a ``large'' field such as the complex numbers.  Beginning with
the work of Bloch [Bl3] and Beauville [Be1,2], the Chow ring of an abelian
variety was found to have some structure arising from a filtration
defined using Pontryagin products or the Fourier transform. In
recent years, several people (Beilinson, Bloch,
Murre, Nori, S.
Saito, H. Saito;
see [B1], [Bl4], [Mu], [N], [S1], [SH]) have
made conjectures asserting the existence of certain
filtrations on the Chow groups of any smooth projective variety.
These conjectures would imply that
the Chow groups have much more structure than has been
hitherto uncovered.  In this paper we define and study a
filtration arising
from the cycle map into $l$-adic \'etale cohomology.  It is
somewhat arithmetic in nature because it is necessary to
first define it over
a field which is finitely generated over its prime subfield.
For varieties
over larger fields, it may be defined by passing to the
limit
 over finitely generated subfields.  Others have studied
this filtration, in
 particular Beauville (unpublished) and Jannsen (see [J4]).
With their permission, I have included some of their results
in this paper.

      Each of the filtrations defined so far has its
advantages and disadvantages.  The conjectural filtration of
Beilinson arising from a spectral sequence
relating extensions of mixed motives and motivic cohomology
would imply many outstanding and deep conjectures about
algebraic cycles.  We will not repeat
this story in detail here since it is treated very clearly
in the paper of Jannsen [J4].  The conjectural filtration of
Murre is known to exist in some cases, and it would also
imply
some of these conjectures, but to establish its existence in
general
appears to require the resolution of some very difficult
conjectures.
The filtrations of H. Saito and S. Saito are
defined in a very explicit way, but it seems difficult to do
specific calculations.  The filtration in this paper has the
advantage that
it may be ``computed'' in some cases, and the powerful tools
of
arithmetic algebraic geometry may be brought to bear on it.
However, it
has two basic
disadvantages.  First, it is not at all obvious that it is independent of
$l$.  Second, suppose $X$ is defined over
a field $k$ which is
finitely generated over its prime subfield and $l$ is a
prime number different from the characteristic of $k$.  Jannsen [J1]
has defined a cycle map into continuous \'etale cohomology (see \S 1):

$$c_n:  CH^n(X)\otimes \bQ\to H^{2n}(X,\bQl(n)).$$

The following is a more general version of a conjecture first made by Soul\'e
over finite fields and by Jannsen and the author over number fields:\\

\noi{\bf Conjecture:}  $c_n$ is injective.\\

This is known when $n$=1 (a consequence of the finite
generation of $CH^1(X)$), but except for special classes of
varieties, we have little information for $n\geq 2$.\\

The Hochschild-Serre spectral sequence:

$$H^r(k,H^q(\Xb,\bQl(n)))\Longrightarrow
H^{r+s}(X,\bQl(n))$$

\noi defines a filtration on $H^{2n}(X,\bQl(n))$, which we
pull back to $CH^n(X)\otimes \bQ$
via the cycle map $c_n$.  We call this the {\it $l$-adic
filtration} on $CH^n(X)\otimes \bQ$.  By definition, the
kernel of the cycle map is contained in all the steps of the
filtration.  For this reason, until we have more information
on this kernel, we can only get results about the filtration
modulo the kernel. We can define another filtration
$\cF^{\bullet}CH^n(X)\otimes \bQ$ on the image of the cycle
map by taking the intersection
 with the filtration obtained from the Hochschild-Serre
spectral sequence.
 It is really this filtration about which we can say
something.\\

          In section 2 of this paper we refine the cycle map
$c_n$ to
get {\it higher $l$-adic Abel-Jacobi maps}
from certain parts of the Chow group to Galois cohomology
groups.  For example,
suppose that $n$=dim $X$ and let $A_0(X)$ denote the group
of zero-cycles of
degree zero modulo rational equivalence.  Let $Alb(X)$ be
the group of $k$-points of the
Albanese variety of $X$ and $T(X)$ the kernel of the
Albanese map:
 $$A_0(X)\to Alb(X).$$
\noi Then we define a map:
$$d_n^2: T(X)\otimes \bQ\to H^2(k,H^{2n-2}(\overline
X,\bQl(n))).$$

\noi This map generalizes the $l$-adic Abel-Jacobi map of
Bloch [Bl1], which in
this case is essentially the map obtained from Kummer theory
on
$Alb(X)$ (see \S 2 for more details). \\

	For the results of this paper we will need to assume
that if $X$ is
a smooth projective variety over a local field $k$ with
ring of integers ${\cal O}$, then $X$  has a regular
proper model over $\cal O$.  Assuming this, the main result
of this paper is:\\

\noindent {\bf Theorem 0.1}:  Let $k$ be a function field in
one variable over a finite field of characteristic $p\neq l$
and
 $X$ a smooth projective variety of dimension $n$  over $k$.  Then the map
$d_n^2$ is zero.\\

\noi Since a global field is of cohomological dimension 2
for
$\bQl$-cohomology, we have:\\

\noi {\bf Corollary}:  For $X$ as above of dimension $n$, we
have
$\cF^iCH^n(X)\otimes \bQ=0$ for all $i\geq 2$.\\

	The following remarks may help clarify the meaning of
these results.
First of all, Beilinson and Bloch have made the conjecture
that if $k$ is a global field and
$X$ is a smooth projective variety over $k$ then the
Albanese kernel $T(X)$ is a
torsion group.  Thus the theorem provides ``evidence'' for
this conjecture for $k$ a
global field of positive characteristic by showing that the
cycle classes of
elements of $T(X)$ in $H^{2n}(X,\bQl(n))$ are trivial.  This
should be true for $k$ a number field as well,
but we can only prove it when $H^{2n-2}(\Xb,\bQ_l(n))$ is
generated by the cohomology classes of
algebraic cycles.  In this case, Beauville showed me an easy
proof
which works over any field.\\

	One of the key technical tools in this paper is the theory of
continuous $l$-adic \'etale cohomology and the Hochschild-Serre
spectral sequence relating
the cohomology of a variety $X$ over a field $k$ with the cohomology
over a separable closure $\kb$ (see \S 1).  Such spectral
sequences do not exist in general for the usual $l$-adic
cohomology defined
by taking the inverse limit of cohomology with finite
coefficients.  Also, we rely on
results of Jannsen on the Galois cohomology of global fields
with
coefficients in $l$-adic cohomology groups of algebraic
varieties over $\overline k$ [J2].  Another key
ingredient is Thomason's purity theorem for $\bQl$-cohomology of schemes of
finite type over the
ring of integers
in a local field of residue characteristic different from
$l$ [T].  This enables us
to construct cycle maps over rings of integers of local
fields.\\

        The $l$-adic Abel-Jacobi maps $d_n^i$ were
discovered independently (and
earlier) by Beauville and by Jannsen, who has also given a
specific description
of the map $d_n^2$ in
terms of extension groups [J3]. \\

At least for $X$ smooth over a field, $CH^*(X)\otimes \bQ\cong
K_0(X)\otimes\bQ,$ and so the $l$-adic filtration
may also be viewed as a filtration on $K_0(X)\otimes \bQ$.  There are Chern
class maps:

$$K_{j}(X)\otimes \bQ\to \bigoplus_{i,j}H^{2i-j}(X,\bQl(i))$$

\noi which one can hope are injective over an absolutely finitely
generated field.  These can be used to define similar filtrations on all
$K$-groups.  We hope to develop this line of
investigation
further in another paper. \\

        I would like to thank U. Jannsen for explaining some
of his unpublished
results to me and allowing me to include some of them in
this paper.  J.-L. Colliot-Th\'el\`ene listened to and helped me with
several points in this paper, especially Proposition 3.2.
A. Beauville informed me that one can easily prove that the
$l$-adic
 filtration is of length at most $n$ on $CH^n$ modulo the
kernel of the
 cycle map.
He also showed me a much simpler proof than the one I
had originally proposed of the triviality of the
higher Abel-Jacobi map in the case where the cohomology of
$\Xb$ is generated by the
classes of algebraic cycles.   I am grateful to him for
allowing me to
include these here.  The referee made several helpful comments which have
improved the exposition.    Finally, I would like to thank
K. Paranjape and S. Saito for useful
discussions.\\

During the long period of thinking about and
writing this
paper, I enjoyed the hospitality of the Max-Planck Institut
f\"ur
Mathematik in Bonn, Universit\"at zu K\"oln, Universit\"at
M\"unster,
Universit\'e de Paris-Sud and the CRM of the Universitat
Aut\`onoma de
Barcelona.  During my one-week visits to the aforementioned
German
universities, I was supported by the Deutsche Forschungsgemeinschaft.

\section{Notation and Preliminaries}
\parindent=0cm

	We denote by $X$ a smooth, projective, geometrically
connected variety over a field $k$ and $l$
 a prime number different from $\mbox{char.} k$.   Let $\kb$ be a
separable closure of $k$
and $G=Gal({\overline k}/k), \Xb=X\times_k\overline k$.  If
$k$ is a global field and $S$ is a
finite set of places of $k$, $G_S$ denotes the Galois group
of a maximal extension of $k$ which is unramified outside
$S$.  For ${\cal X}$ a noetherian scheme, we denote by $CH^n({\cal X})$
 the Chow group of codimension $n$ cycles modulo rational
equivalence.  For $X$ proper and geometrically connected over a field,
$CH_0(X)$ is the group of zero-cycles modulo
rational
equivalence and $A_0(X)$ the group of zero-cycles of degree
zero modulo
rational equivalence.  The group of $k$-points of the
Albanese variety
of $X$ will be denoted by $Alb(X)$ and $\alpha : A_0(X)\to
Alb(X)$
denotes the Albanese map.  We
denote the kernel of this map by $T(X)$.

If $M$ is an abelian group upon which $G$ acts continuously,
we denote
the continuous Galois cohomology groups of $G$ with values in $M$ by
$H^*(k,M)$ (see \S 1 for a definition of these groups). The subgroup of $M$
consisting of elements killed by $l^m$
is denoted by $M_{l^m}$ and $$T_lM={\plim}M_{l^m}.$$   We set
$V_lM=T_lM\otimes_{\bZl}\bQl.$  As usual, we denote by $\bZl(1)$ the
$l$-adic sheaf of $l$-primary roots of unity.  For $n>0, \bZl(n)$
denotes $\bZl(1)^{\otimes n}$; for $n<0, \bZl(n)=Hom(\bZl(-n),\bZl)$.
 We use similar notation for $\bQl$-sheaves.\\

	Let $X$  be a smooth variety over a field $k$ and $l$ a
prime number different from the characteristic of $k$.
Jannsen has defined a cycle map
[J1]:
$$CH^n(X)\otimes \bQ \to H^{2n}(X,\bQl(n)),$$

\noi where the group on the right is the continuous $l$-adic
cohomology group with
values in the $l$-adic sheaf $\bQl(n)$.  We briefly recall the
definition of these groups.
Let $({\cal F}_m)$ be a projective system of sheaves of
$\bZ/l^m$-modules
on $X$.  The
functor
which sends $(\cF_m)$ to $$\plim H^0(X,(\cF_m))$$ is left
exact, and we define
$H^i(X,(\cF_m))$ to be the $i$-th right derived functor.
When the inverse system is $\bZ/l^m(j)$,
we denote these groups by $H^i(X,\bZl(j))$.  Continuous
\'etale cohomology with
$\bQl$-coefficients is defined by tensoring the $\bZl$-cohomology with $\bQl$.

 We shall need the following properties of
this cohomology which can be found in Jannsen's paper
[J1].\\

\noindent (i)  There is an exact sequence:

$$0\to {\plim} ^1H^{i-1}(X,(\cF_m))\to H^i(X,(\cF_m))\to
\plim H^i(X,\cF_m)\to 0.$$

In particular, if $H^i(X,\bZ/l^m(j))$ is finite for all $i$
and $m$ then
$$H^i(X,\bZl(j))\to \plim H^i(X,\bZ/l^m(j))$$
is an isomorphism.  This is true if, for example, $k$ is
separably closed, finite or local.  This is not true in
general if $k$ is a number field.\\

\noindent (ii) Let $\overline k$ be a separable closure of
$k$, $G=Gal(\overline k/k)$
and $\overline X=X\times _k\overline k$.  Then there is a
Hochschild-Serre
spectral sequence:
$$E_2^{p,q}=H^p(G,H^q(\overline X,\bZl(j)))\Longrightarrow
H^{p+q}(X,\bZl(j)).$$

\noindent {\bf Theorem 1.1} (Jannsen):  Let $X$ be a smooth
projective variety over a field $k$.
Then the Hochschild-Serre spectral sequence with $\bQl$-coefficients
degenerates at $E_2$.\\

	Since this result has not been published, let us
indicate the ingredients in the proof.  The idea is to use
Deligne's criterion for a spectral sequence to degenerate
[D1].  In order to use this in this situation, one needs a
good
derived category of $\bQl$-sheaves.  One cannot use
Deligne's
definition given in ([D2], 1.1.2d) because the cohomology groups with
finite
coefficients need not be finite.  Jannsen has given a good
definition
of such a derived category in the general case.  Once this
is done one applies the Hard Lefschetz
Theorem to satisfy Deligne's criterion.  Ekedahl has
published another
definition of the derived category of $\bQl$-sheaves [E], and
this should
also suffice to prove the degeneration of the Hochschild-Serre spectral
sequence.

\section{The $l$-adic filtration}

{}From the Hochschild-Serre spectral sequence for continuous
cohomology:
$$E_2^{p,q}=H^p(k,H^q(\overline X,\bQl(n)))\Longrightarrow
H^{p+q}(X,\bQl(n)),$$
we get a filtration on $H^{2n}(X,\bQl(n))$:
$$H^{2n}(X,\bQl(n))=F^0\supset F^1\supset . . .\supset
F^{2n}$$
and $F^i/F^{i+1}=E_\infty ^{i,2n-i}$.\\

{\bf Definition 2.1}: Let $X$ be a smooth, projective,
geometrically connected variety over a field $k$ which is
finitely generated over its prime subfield.  Then we define:

$$F^iCH^n(X)\otimes \bQ=c_n^{-1}[F^iH^{2n}(X,\bQ_l(n))],$$
where $c_n$ is the cycle map (see \S 1).\\

If $k$ is not finitely generated, we define:

$$F^iCH^n(X)\otimes
\bQ=\lim_{\stackrel{\rightarrow}{L}}F^iCH^n(X_L)\otimes \bQ,$$

where $L$ runs over all finitely generated fields contained in $k$ for which
$X$ has a model over $L$.\\

Finally, set $$\cF^iCH^n(X)\otimes \bQ=\mbox{Image}\,c_n \cap
F^{i}H^{2n}(X,\bQ_l(n)).$$  Of course, {\it a priori}
$\cF^{\bullet}CH^n(X)\otimes\bQ$ is only a filtration on the
image of the cycle
map.\\

Let $X$ be defined over a field $L$ which is finitely generated over the prime
subfield.  By the definition of the filtration on $CH^n(X)\otimes \bQ$
and the
degeneration of the Hochschild-Serre spectral sequence (Theorem 1.1), we get
maps:

$$d_n^i:\, F^iCH^n(X)\otimes \bQ\to H^i(k,H^{2n-
i}(\Xb,\bQl(n)))$$

which we call the {\it higher $l$-adic Abel-Jacobi
mappings}.  For $X$ over arbitrary $k$, we define the higher Abel-Jacobi map as
the limit of the maps defined above for $X_L$, where $L$ runs over finitely
generated subfields of $k$ over which $X$ has a model. \\

When
$i$=1, these are the $l$-adic Abel-Jacobi mappings defined
by Bloch
[Bl1].  When $n=\mbox{dim}X$, this is the Albanese map
followed by the
map given by Kummer theory for $Alb(X)$, as we now
explain.\\

Let $A_0(X)$ be the subgroup of $CH^n(X)$ consisting of
cycles of degree zero
and let $Alb(X)$ denote the group of $k$-points of the
Albanese variety of $X$.  Consider the Albanese map $$\alpha:  A_0(X)\to
Alb(X)$$  and its kernel $T(X)$.  Then
we have a commutative diagram with exact rows:
$$\matrix{0 &\to & F^2&\to & F^1&\to& H^1(k,H^{2n-
1}(X,\bQ(n)))&\to & 0\cr
&&f\uparrow&&\uparrow&&\uparrow\cr
0&\to& T(X)\otimes \bQ&\to& A_0(X)\otimes \bQ&\to
&Alb(X)\otimes \bQ&\to &0}.$$
Here the right vertical map is given as follows.  By Poincar\'e duality and the
Weil pairing, we have an isomorphism:

$$H^{2n-1}(\overline X,\bQl(n))\cong V_l(Alb(\overline X)),$$

where $V_l$ denotes the Tate $\bQl$-vector space of $Alb(\overline X)$.  Then
the right vertical map is defined
by doing Kummer theory modulo $l^m$ on
$Alb(X)$, passing to the inverse limit over $m$ and tensoring with $\bQl$.  The
commutativity of the right square is proved in the
appendix to this paper.  This induces the map $f$.
Composing $f$ with the map obtained from the spectral
sequence:
$$F^2\to E_\infty ^{2,2n-2},$$
we get a map
$$d_n^2:\, T(X)\otimes \bQ\to H^2(k,H^{2n-2}(\overline
X,\bQl(n))).$$

{\bf Remarks 2.1.1}:  (i)  It would be interesting to find a complex
analytic analogue of the maps $d_n^i$.  Of course, $d_n^1$ may be
defined in the complex category as the Abel-Jacobi map of Griffiths.
The target group of our $d_n^i$ may be reinterpreted as
$Ext_{\bQl}^i(\bQl,H^{2n-i}(\Xb,\bQl(n)))$ in the category of continuous
Galois modules.  The complex analogue of this is
$Ext_{MHS}^i(\bQ,H^{2n-i}(X,\bQ))$, where $MHS$ denotes the
category of mixed Hodge structures.  For $i=1$, this works fine, as it
is well-known that $Ext_{MHS}^1(\bZ,H^{2n-1}(X,\bZ))$ for $X/{\bf C}$ is
the intermediate Jacobian of Griffiths. Unfortunately, all the higher
$Ext$-groups are zero.  Several people have suggested
taking extensions in a suitable category of variations of Hodge
structure.\\

(ii)  The filtration $F_B$ of S. Saito (see [Sa3], Definition
(1-4)) is contained in the $l$-adic filtration.  It is not
clear to me whether this is the case for the filtration $F_{BM}$ defined in
(loc. cit.
Definition (1-5)).  For much
more on the relations between various filtrations, see the paper of
Jannsen [J4].\\

(iii)        The $l$-adic Abel-Jacobi maps defined above have a conjectural
mixed motivic interpretation.  Recall that
there are motivic cohomology groups with $\bQ$-coefficients
$H_{\cal M}^i(X,\bQ(j))$ with
$$H_{\cal M}^{2n}(X,\bQ(n))=CH^n(X)\otimes \bQ.$$
These may be defined by:
$$H_{\cal M}^i(X,\bQ(n))=gr_{\gamma}^{n}K_{2n-i}(X)\otimes \bQ,$$
where $\gamma$ denotes the gamma filtration in algebraic
$K$-theory
 (here it
is crucial that we are dealing with {\bf Q}-coefficients
because the integral
theory is much less developed).  Beilinson [B1] has made the
conjecture that
these groups have a canonical filtration coming from a
spectral sequence:
$$Ext^p_{\cal M\cal M}(\bQ,H^q(X,\bQ(j)))\Longrightarrow
H_{\cal M}^{p+q}(X,\bQ(j)).$$
Here the subscript $\cal M\cal M$ denotes the conjectural
category of mixed motives on $X$ and the $Ext$ is taken in
this category.
The maps given above are in some sense ``$l$-adic
realizations'' of the corresponding maps obtained
from this conjectural spectral sequence.  This point
of view is very useful for learning what to expect.  In fact,
if Beilinson's
conjecture on the existence of filtrations on Chow groups
([B1], [J4]) is true,
then the filtration defined in this paper must agree with
his (see [J4],
Lemma 2.7).\\

We now prove some basic results on these filtrations.\\

{\bf Proposition 2.2} (Beauville):  For $X$ as above, we
have

$$F^{n+j}CH^n(X)\otimes \bQ=F^{n+1}CH^n(X)\otimes \bQ$$

for all $j\geq 1.$ \\

{\bf Proof}: It suffices to
show that the higher Abel-Jacobi maps $d_n^i$ are zero for $i>n$.
To prove this, let $d$ be the dimension of
$X$ and let $Y$ be a smooth hyperplane section of $X$.  Then
the diagram:

$$\matrix{F^{n+j}CH^n(X)\otimes \bQ&\to&H^{n+j}(k,H^{n-
j}(\Xb,\bQl(n)))\cr
\downarrow&&\downarrow\cr
F^{n+j} CH^{d+j}(X)\otimes \bQ&\to&H^{n+j}(k,H^{2d-
n+j}(\Xb,\bQl(d+j)))}$$

is commutative. Here the left vertical arrow is given by
intersecting
cycles on $X$ with the class of the $(d+j-n)$-fold
intersection of $Y$
with itself, and the right vertical arrow is given by
cupping
$(d+j-n)$-times with the class of $Y$ in $H^2(X,\bQl(1))$.
This
commutativity is easily seen from the compatibility of the
Hochschild-Serre spectral sequence and the cycle map with
intersection
products ([J1], Lemma 6.14).  Now since $j\geq 1$, the
bottom left group is zero.  But the right vertical arrow is
an isomorphism by the
Hard Lefschetz Theorem, and hence the image of the top
horizontal arrow is trivial.  This proves the propostion.\\

Recall that the cohomology group $H^i(\Xb,\bQl)$ is said to
be {\it
algebraic} if it is generated by Tate twists of cohomology
classes of
algebraic cycles.  Thus if $i$ is odd, this means that the
group is
zero, and if $i=2j$ is even then the cycle map:

$$CH^j(\Xb)\otimes \bQl\to H^{2j}(\Xb,\bQl(j))$$

is surjective (we keep the Tate twists because we will need
them later).
A basic conjecture of Bloch ([Bl2], Conjecture 0.4 in the
case of
surfaces) predicts that if $X$ is defined over
$\bf C$ and the cohomology of $X$ is supported on codimension 1 subschemes,
then the
Albanese map:

$$\alpha:\,  A_0(X)\to Alb(X)$$

is an isomorphism.  We now give an argument of Beauville to
show how a weaker form of
this conjecture is implied by the conjecture on injectivity
of the
cycle map for varieties over finitely generated fields (see the introduction to
this paper).\\

{\bf Proposition 2.3} (Beauville):  Let $X$ be a smooth
projective
variety of dimension $d$ over a field $k$.  Let $i$ be a
positive integer and assume that
the cohomology groups $H^{2n-2i}(\Xb,\bQl(i))$ and
$H^{2d-2n+2i}(\Xb,\bQl(d-n+i))$ are
generated by algebraic cycles.  Then the map:

$$F^{2i}CH^n(X)\otimes\bQ\to H^{2i}(k,H^{2n-2i}(\Xb,\bQl(n)))$$

is zero. \\

{\bf Proof}:   First note that for any finite extension $L/k$, the natural
map:

$$H^{2i}(k,H^{2n-2i}(\Xb,\bQl(n)))\to H^{2i}(L,H^{2n-2i}(\Xb,\bQl(n)))$$

is injective because its kernel is killed by $[L:k]$.  By the algebraicity
assumption, there is a finite extension $L/k$ over which
$G=Gal(\kb/L)$
acts trivially on $H^{2n-2i}(\Xb,\bQl(n-i))$ and on
$H^{2d-2n+2i}(\Xb,\bQl(d-n+i))$.  Thus we may and do assume that $Gal(\kb/k)$
acts
trivially on these groups.
Let $Y_1,...Y_{\nu}$ be codimension $d-n+i$
cycles on $X$ whose cohomology classes
$\xi_1,...,\xi_\nu$ form a basis for $H^{2d-2n+2i}(\Xb,\bQl(d-n+i))$.  Then
the following
diagram is commutative:

$$\matrix{F^{2i}CH^n(X)\otimes \bQ&\to&H^{2i}(k,H^{2n-2i}(\Xb,\bQl(n)))\cr
\downarrow&&\downarrow\cr
[F^{2i}CH^{d+i}(X)\otimes\bQ]^{\nu}&\to&[H^{2i}(k,H^{2d}(\Xb,\bQl(d+i)
))]^{\nu}}.$$

Here the left vertical arrow is given by intersecting with
the class of
$Y_j$ on the $j$-th factor.  The right vertical arrow is
given by cupping
with the class of $\xi_j$ on the $j$-th factor; it is an
isomorphism by the assumption that the
cohomology is algebraic
and Galois acts trivially on $H^{2n-2i}(\Xb,\bQl(n-i))$.
The
commutativity is easily established by the compatibility
of the Hochschild-Serre spectral sequence and the cycle map
with
intersection product.  Now the bottom left group is zero if
$i\geq 1$.  Hence the top
horizontal map is zero.  This completes the proof of the
proposition. \\

{\bf Corollary 2.4}:  Let $X$ be a smooth projective variety of dimension
$d$ over an algebraically closed field $k$.
assume that for $d\leq i\leq 2d-2$,
the cohomology group $H^i(\Xb,\bQl)$ is algebraic.  Let $K$ be a
finitely generated field of definition of $X$ and assume that the
cycle map:

$$c_d:\, CH^d(X_M)\otimes \bQ\to H^{2d}(X_M,\bQl(d))$$

is injective over any finitely generated field $M$ containing $K$.    Then the
Albanese map:

$$\alpha:\: A_0(X)\to Alb(X)$$

is an isomorphism.  In particular, if $k={\bf C}$ and $X$ is a surface with
$H^0(X,\Omega^2_X)=0$, then injectivity of the cycle map $c_2$ over all
finitely generated fields
containing $K$ implies Bloch's conjecture.\\

{\bf Proof}:  Let $a$ be in the kernel of $\alpha$.  Then
$a$ is defined
over some finitely generated field $M$.  The assumption on
injectivity of the cycle map implies
that the maps:

$$F^iCH^d(X_M)\otimes \bQ\to H^i(M,H^{2d-i}(\Xb,\bQl(d)))$$

are injective.  If $i\geq 2$ then the assumption that the
cohomology is
algebraic and Proposition 2.3 imply that the map is zero.
When $i=1$,
the map may be identified with $\alpha$ over $M$ followed by Kummer
theory for $Alb(X_M)$ (see the appendix to this paper).
Hence $\alpha$ is injective up to torsion.  But by Roitman's
theorem
[Ro], [M2],
the kernel of $\alpha$ has no torsion, so it is zero.  The last statement
follows immediately from the exponential sequence, Serre duality and the
comparison theorem between singular and \'etale cohomology.\\

{\bf Remark 2.4.1}:  This result should be true under the weaker hypothesis
that for each $d\leq i\leq 2d-2$, the cohomology group $H^i(\Xb,\bQl)$ is
supported on codimension one subschemes.  This is Bloch's conjecture in the
general case.\\

{\bf Proposition 2.5}:  Let $X$ be a smooth, projective
variety over a
field $k$ which is finitely generated over its prime
subfield and $l$
a prime number different from the characteristic of $k$.
Then the image
of the cycle map:

$$c_n:  CH^n(X)\otimes \bQl\to H^{2n}(X,\bQl(n))$$

is a finite dimensional $\bQl$-vector space.\\

{\bf Proof}:   There exist an open subset $U$ of $\mbox{Spec}\bZ[1/l]$, a
regular
ring $A$ of finite type and smooth
over $U$ with fraction field $k$, and a smooth proper
scheme $\sX$
over $A$ with generic fibre $X$.  There is a cycle map:

$$CH^n(\sX)\otimes \bQl\to H^{2n}(\sX,\bQl(n))$$

(see e.g. ([S2], Lemma 5.3) or Proposition 2.9 below).  Now
$H^{2n}(\sX,\bQl(n))$ is a finite dimensional
$\bQl$ vector space (see e.g. [CT-R], Proof of Th\'eor\`eme
1.1, (d)) and since the map

$$CH^n(\sX)\otimes \bQl\to CH^n(X)\otimes\bQl$$

is surjective, we get the proposition.\\

{\bf Corollary 2.6}:  With hypotheses as in Proposition 2.5,
for any
nonnegative integer $i$, the image
of the higher $l$-adic Abel-Jacobi map:

$$F^iCH^n(X)\otimes\bQl\to H^i(k,H^{2n-i}(\Xb,\bQl(n)))$$

is a finite dimensional $\bQl$ vector space.\\

Proposition 2.5 and Corollary 2.6 are not very refined
results since
their proofs only use the existence of a smooth proper model
of $X$ over
an open subset of a model of $k$ over an open subset of
Spec$\bZ$.  In order to get finer
results, we need to assume the existence of a regular proper
model
$\sY$ of $k$ over $\bZ$ together with a regular proper model
of $X$
over $\sY$.  In the rest of this section we assume we have
such
models.  Our philosophy is that the image of the higher
Abel-Jacobi map
$d_n^i$ lies in the {\it unramified part} of
$H^i(k,H^{2n-i}(\Xb,\bQl(n))),$ in the following sense.\\

{\bf Definition 2.7}:  Let $k$ be a field which is finitely generated over its
prime subfield, $l$ a prime number which is invertible in $k$ and ${\cal Y}$ a
regular proper model of $k$ over $\bZ$.  Let $F$ be an $l$-adic sheaf on ${\cal
Y}$. We denote by $H^i_{nr}(k,F)^{(l)}$ the $\bQl$-vector space:

$$\mbox{Ker}[H^i(k,F)\to \prod_{y\in {\cal Y}^{1\prime}}H^{i+1}_y(A_y,F)].$$

Here ${\cal Y}^{1\prime}$ is the set of codimension one points of ${\cal Y}$
which do not lie above $l$ and $A_y$ is the valuation ring of $y$.  By abuse of
notation, we also denote by $F$ the sheaves on $k$ and $A_y$ induced by $F$.
If $F$ is a locally constant sheaf and $l$ is invertible on ${\cal Y}$, it
follows easily from purity that this is the usual definition of unramified
cohomology (see [CT-O], for example).\\

Now let $X$ be a smooth projective variety over $k$.  Assume
that we have
a regular proper scheme $\sY$ over $\bZ[1/l]$ with function
field $k$ and a
regular proper scheme $\sX$ together with a morphism:

$$f:\,\sX\to \sY$$

whose generic fibre is $X/k$.  We consider the
 $l$-adic
sheaf
$F=R^{2n-i}f_*\bQl(n)$, which is locally constant on the smooth locus of $f$,
and we can
define the
unramified cohomology of $k$ with values in this sheaf.
Then we
have:\\

{\bf Conjecture 2.8}:  With notation as above, the image of
the higher
Abel-Jacobi map $d_n^i$ is contained in $H^i_{nr}(k,F)^{(l)}$.  \\

{\bf Remark 2.8.1}:  In the definition of unramified cohomology and the
conjecture, we have ignored
points of $\sY$ lying above $l$,
not because we don't think these are important, but rather
because
we are unable to formulate an intelligent conjecture on the
image of the
Abel-Jacobi map at these places.  What is needed is a higher
dimensional
analogue of the finite part, exponential part and geometric
part of
Bloch-Kato [BK].  This will involve $l$-adic cohomology over
a complete
discretely valued field with nonperfect residue field of characteristic $l$,
and
this is not
very well developed at the moment.\\

 Another way to formulate Conjecture 2.8 is that the
filtration on the
Chow groups defined in 2.1 should extend to $\sX$.  We define
this
filtration using the following result:\\

{\bf Proposition 2.9}:  Let $R$ be a ring which is of type (i), (ii) or
(iii) below:\\

(i)  of finite type over $\bZ$\\

(ii)  a localization of a ring as in (i)\\

(iii)  the henselization of a ring as in (ii)\\

Let $\sX$ be a regular scheme of
finite type over $R$ and assume that $l$ is
invertible on $\sX$. Then there is a cycle
map:

$$c_n: CH^n(\sX)\to H^{2n}(\sX,\bQl(n)).$$

{\bf Proof}:  Thomason [T] has proved a purity theorem for
$\bQl$-cohomology
of such schemes, and then the existence of a cycle map may
be proven by
copying the proof in e.g. ([M1], Chapter VI, \S 6).\\

{\bf Definition 2.10}:  Let $f:\,\sX\to \sY$ be a proper
morphism of schemes of
finite type over a ring $R$ as in 2.9 and $l$ a prime number invertible on
$\sY$.  Then
we have the Leray spectral sequence for continuous
cohomology ([J1],
3.10):

$$H^p(\sY,R^qf_*\bQl(n))\Longrightarrow
H^{p+q}(\sX,\bQl(n)).$$

This defines a filtration on $H^{2n}(\sX,\bQl(n))$ which we
can pull
back to\newline $CH^n(\sX)\otimes \bQ$ via $c_n$.  We denote the
resulting filtration
$F^{\bullet}CH^n(\sX)\otimes \bQ$.  Of course, this filtration may
depend on $\sY$, but this will be fixed in what follows.\\

The following is a properness conjecture for the filtration
on a regular
proper scheme over a discrete valuation ring.\\

{\bf Conjecture 2.11}:  Let $A$ be a discrete valuation ring as in 2.9(ii) or
(iii)
with
fraction field $k$ and $l$
a prime number invertible on $A$.  Let $\sX$
be a regular proper scheme of finite type over $A$ with
generic fibre
$X$.  Then the natural map:

$$F^iCH^n(\sX)\otimes \bQ\to F^iCH^n(X)\otimes \bQ$$

is surjective.\\

{\bf Proposition 2.12}:  With notation as in 2.11, assume
that $A$ is
Henselian and $\sX$ is
smooth over $A$.  Then Conjecture 2.11 is true for $\sX$.\\

{\bf Proof}: We need a lemma.\\

 {\bf Lemma 2.13}:  Let $A$ be a Henselian discrete
valuation ring with
 fraction field $k$ and $(F_m)$ a projective system of
locally constant $l$-primary
 sheaves on $A$.  Then the natural map on continuous
cohomology groups:

$$H^i(A,(F_m))\to H^i(k,(F_m))$$

is injective.\\

{\bf Proof}:  Let $s$ be the closed point of $S=\mbox{Spec}A$.  For each $m$ we
have an
exact sequence:

$$(*)...\to
H^{i-1}(k,F_m)\stackrel{f}{\to} H^{i}_{s}(A,F_m)\to
H^i(A,F_m)\to H^i(k,F_m)\to ...$$

Now the second group in this sequence is isomorphic to
$H^{i-2}(s,F_m(-1))$ by purity for locally constant sheaves
on $A$.
The map $f$ is split surjective with a splitting being given
by the map
($\pi$ a uniformizing parameter of $k$):

$$H^{i-2}(s,F_m(-1))\cong
H^{i-2}(A,F_m(-1))\stackrel{\cup\pi}{\to}H^{i-1}(k,F_m).$$

Hence the map obtained from $f$ above by passing to the
inverse limit
over $m$ is surjective as well.  A similar statement holds
for the case
of the map
$${\plim}^1H^{i-1}(k,F_m)\to {\plim}^1H^{i-2}(s,F_m(-1)).$$
Then by the property (i) of continuous cohomology discussed
in \S 1
above, the
map:

$$H^{i-1}(k,(F_m))\to H^{i-2}(s,(F_m))$$

is surjective.  Now there is a localization sequence like
(*) above for
continuous cohomology with support and a purity theorem for
the
cohomology with support in $s$ ([J1], 3.8
and same argument as in Theorem 3.17 of that paper).  We conclude that the map

$$H^i(A,(F_m))\to H^i(k,(F_m))$$

is injective.  This completes the proof of the lemma.\\

Returning to the proof of the proposition,
since $\sX$ is smooth over $A$, the sheaves
$F_{i,m}:=R^if_*\bZ/l^m(n)$ are locally constant on $A$
(smooth and
proper base change theorem), so the maps

$$H^j(A,(F_{i,m}))\to H^j(k,(F_{i,m}))$$

are injective.  Then the Leray spectral sequence
$$H^p(A,R^qf_*(F_m))\Longrightarrow
H^{p+q}(\sX,(F_m))$$

degenerates at $E_2$ because it does so for $X/k$
(Theorem 1.1).  Now an easy induction starting from the fact
that the
 map

 $$CH^n(\sX)\to CH^n(X)$$

is surjective proves the proposition.\\

{\bf Remarks 2.14}:  (i)  When $\sX$ is not smooth over $A$,
it seems much more difficult to prove
properness of the filtration.  In fact, it is not even known
if the map
$$F^1CH^n(\sX)\otimes \bQ\to F^1CH^n(X)\otimes \bQ$$

is surjective.  In Proposition 3.2 below we prove this when
$n=\mbox{dim}X$,
$A$ is Henselian and $k$ with finite or separably closed residue field.\\

(ii)  It should be possible to prove Proposition 2.12 for
any discrete
valuation ring $A$ using results of Gillet [Gi].  His method
may be used to
show that if $F$ is a locally constant sheaf on $A$ with
finite stalks, then the natural
map:

$$H^i(A,F)\to H^i(k,F)$$

is injective.  We need this for continuous cohomology.   In some sense,
our result for Henselian rings is not what we really want because the
fraction field $k$ is usually not finitely generated over the prime
field.  However, as
will be
shown in \S 3, in practice the Henselian case is very
helpful.\\

\section{Reductions in the proof of the Main Theorem}
        In this section we show how to reduce Theorem 0.1 to
a Hasse
principle for Galois cohomology which has been proved in
some cases by
Jannsen. In order to do this, we need some local
considerations and,
unfortunately, a hypothesis:\\

{\bf Hypothesis 3.1}:  Let $A$ be a discrete valuation ring
with fraction field $k$ and $X$
a smooth projective variety over $k$.  \\

  (Resolution of singularities) $X$ has a regular
proper model $\sX$ over $A$. \\

Hypothesis 3.1 will be needed in order to have a cycle map
for $\sX$
(see Proposition 2.9) and for Proposition 3.2 below.\\

The following result may be found in ([B1], Lemma 2.2.6(b)).
We give some
of the details of the proof outlined there as well as the
slight
generalization we require.\\

{\bf  Proposition 3.2}:  Let $A$ be a Henselian discrete
valuation ring with fraction field $k$ and  residue field
$\bF$ which is of
characteristic different from $l$. Assume that $\bF$ is either finite or
separably closed.  Let $X$ be a smooth,
projective, geometrically connected variety of dimension $n$
over $k$ and $\xi$  an element of
$CH^n(X)\otimes \bQ$ of degree 0.
Then there exists a finite extension $L/k$ with valuation
ring $B$, a
regular proper model $\sX$ of $X_L$ over $B$ and a cycle in
$CH^n(\sX)\otimes \bQl$ which restricts to $\xi_L$ in $CH^n(X_L)\otimes\bQl$
and whose
cohomology class in $H^{2n}(\sX,\bQl(n))$ is zero.\\

{\bf Proof}:  If $\xi$ is a zero-cycle of degree zero on
$X$, its cohomology class is trivial in
$H^{2n}(\Xb,\bQl(n)).$  We claim that its cohomology class
is zero in $H^{2n}(X,\bQl(n))$.
To see this, note that by Lemma 3.3 below, the group $A_0(X)$ is an extension
of a torsion group by a group which is divisible prime to the residue
charactersitic of $k$.  Now the group $H^{2n}(X,\bZl(n))$ contains no
nontrivial $l$-divisible subgroup ([J1], Cor. 4.9).  Thus the restriction of
the integral cycle map to cycles of degree zero:

$$A_0(X)\to H^{2n}(X,\bZl(n))$$

has torsion image, so the class of $\xi$ in $H^{2n}(X,\bQl(n))$ is trivial, as
claimed.\\

Assume for the moment
that $X$ is a
curve. Let $\sX$ be a regular proper
model of $X$ over $A$ with special fibre $Y$.    There is a commutative diagram
with exact top row:

$$\matrix{CH^0(Y)\otimes\bQl&\to&CH^1(\sX)\otimes\bQl&\to&CH^1(X)\otimes\bQl&\to&0\cr
&&\downarrow&&\downarrow\cr
&&H^2(\sX,\bQl(1))&\to &H^2(X,\bQl(1))}.$$

The vertical maps are cycle maps and the top sequence is the exact sequence
for the Chow groups of a scheme, an open subscheme and its complement.
By Lemma 3.4 below, we can
identify $H^2(\sX,\bQl(1))$ with\linebreak $(CH^0(Y_{{\overline \bF}})\otimes
\bQl)^{*\,Gal({\overline \bF}/\bF)}$, and the cycle map

$$CH^1(\sX)\otimes \bQ\to H^2(\sX,\bQl(1))$$

constructed in Lemma 2.9 is given by intersecting with components of the
geometric special fibre.
Here * denotes $\bQl$-linear dual.
Now the intersection pairing on $CH^0(Y)/[Y]$ in
$\sX$ is nondegenerate,
where $[Y]$ is the class of the whole special fibre  ([L],
Ch. 12, \S
4), and it induces an isomorphism:

$$CH^0(Y)/[Y]\otimes \bQl\to [CH^0(Y_{{\overline \bF}})^0\otimes
\bQl]^{*\,Gal({\overline \bF}/\bF)},$$

where the superscript 0 means the kernel of the degree map.

If $\xi$ is a zero-cycle of degree zero on $X$ and
$\tilde{\xi}$ is any lifting
of $\xi$ to a codimension 1 cycle on $\sX$, then the class
of $\tilde{\xi}$ in
$CH^0(Y)$ is of total degree zero, where
total degree means the sum of the degrees on each
irreducible component of $Y$.
Hence we can modify $\tilde{\xi}$ by an element
of $CH^0(Y)\otimes \bQl$ so that its class in\linebreak
$[CH^0(Y_{\overline{\bF}})\otimes
\bQl]^*$ is trivial.
This proves the lemma in the case of a curve.  For the
higher dimensional case,
let $\xi$ be an element of $CH^n(X)\otimes \bQ$ of degree
zero.  Then passing to
 a finite extension of $k$ if necessary (if $k$ is not
perfect), we may find a
 smooth, projective, geometrically connected curve $C$
passing through
 $\xi$ ([C], [AK]).  Let $B$ be the integral closure of $A$
in $L$ and
 let $\cC$ be a regular proper model of $C$ over $B$.  Then
by blowing up, the
 rational map $\cC\to \sX$ may be resolved to get a
morphism:

$$f: \tilde{\cC}\to \sX,$$

with $\tilde{\cC}$ regular.

Let $Z, Y$ be the special fibres of $\tilde{\cC}$ and $\sX$,
respectively.  Then
the following diagram commutes:

$$\matrix{CH^1(\tilde{\cC})\otimes\bQl&\to& CH^n(\sX)\otimes\bQl\cr
\downarrow&&\downarrow\cr
[CH^0(Z_{\overline{\bF}})\otimes \bQl]^*&\to&[CH^0(Y_{\overline{\bF}})\otimes
\bQl]^*}.$$

Here the top horizontal map is the push-forward map and the bottom horizontal
map is defined as the dual of the pullback map on irreducible components.

By the first part of the proof, we can find a cycle
$\tilde{\xi}$ in
$CH^1(\tilde{\cC})\otimes \bQl$ whose image in $CH^1(C)\otimes\bQl$ is
equal to $\xi$
and whose class is trivial in $CH^0(Z)\otimes \bQl$.  Pushing
$\tilde{\xi}$ to
$CH^n(\sX)\otimes \bQl$, we get our desired cycle.  This
completes the proof of
the proposition. \\

We now prove the two lemmas which were used in the proof of Proposition 3.2.\\

{\bf Lemma 3.3}:  Let $k$ be a Henselian discrete valuation field.  Assume that
the residue field $\bF$  is finite or separably closed.
Let $X$ be a smooth projective variety over $k$.  Then the group $A_0(X)$ is
an extension of a torsion group by a group which is divisible prime to the
characteristic of $\bF$.\\

{\bf Proof}:  As in the latter part of the proof of Proposition 3.2, given any
zero-cycle on $X$, after passing to a finite extension of $k$ if necessary, we
may find a smooth, projective, geometrically connected curve going through its
support.  By an easy norm argument, we see that it suffices to prove the result
for the case of curves.  This is implied by the following fact applied to the
Jacobian of the curve:\\

{\bf Fact}:  Let $k$ be as in the hypothesis of Lemma 3.3 and $A$ an abelian
variety over $k$.  Then the group $A(k)$ is an extension of a finite group by
a group which is divisible prime to the characteristic of $\bF$.  \\

This is well-known if the residue field is finite.  To see it when the
residue field is separably closed, let ${\cal A}$ be the N\'eron model of $A$
over the valuation ring $R$ of $k$ with special fibre $\widetilde{{\cal A}}$
(see [BLR], Ch. 1, \S 1.3, Cor. 2 for the existence of ${\cal A}$).  Let
$A_0(k)$ be the subgroup of $A(k)$ consisting of points which specialize to
elements in the connected component of identity  $\widetilde{{\cal A}}_0$ of
$\widetilde{{\cal A}}$.  Then the quotient $A(k)/A_0(k)$ is finite.  The group
$\widetilde{{\cal A}}_0(\bF)$ is the group of points of a connected algebraic
group over the separably closed field $\bF$, and hence is divisible prime to
the characteristic of $\bF$.  Let $A_1(k)$ be the kernel of the surjective
reduction map:

$$A_0(k)\to \widetilde{{\cal A}}_0(\bF).$$

Then $A_1(k)$ is a pro-$p$ group , where $p$ is the characteristic of $\bF$.
Hence it is uniquely divisible prime to $p$.  Putting everything together, we
get the fact.\\

{\bf Lemma 3.4}:  Let $A$ be a Henselian local ring with residue
field $\bF$ which is finite or separably closed.  Let $\sX$ a regular proper
flat integral scheme of dimension $n$ over $A$ with special fibre $Y$. Let $l$
be a prime number different from the characteristic of $\bF$.  Then

$$(\dagger)\:\: H^{2n}(\sX,\bQl(n))\cong [CH^0(Y_{{\overline \bF}})\otimes
\bQl]^{*\,Gal({\overline \bF}/\bF)},$$
where * denotes $\bQl$-dual.
The map

$$CH^n(\sX)\to [CH^0(Y_{{\overline \bF}})\otimes \bQl]^*$$

obtained by composing this isomorphism with the cycle map in
Proposition 2.9 is given by intersecting with components of the
geometric special fibre.\\

{\bf Proof}:  This is known, but I could not find
a good
reference.  By the proper base change
theorem, we have:

$$H^{2n}(\sX,\bQ_l(n))=H^{2n}(Y,\bQl(n)).$$
 This last group may be easily calculated as follows.
Since the \'etale cohomology with $\bQl$-coefficients of $Y$ and $Y_{red}$ is
the same, we may assume that $Y$ is reduced. Let
$Y_{sing}$ be
the singular locus of $Y$.  Then there is an exact sequence
for
cohomology with compact support:

$$...\to H_c^{2n-1}(Y_{sing},\bQl(n))\to H_c^{2n}(Y-
Y_{sing},\bQl(n))\to
H_c^{2n}(Y,\bQl(n))\to H_c^{2n}(Y_{sing},\bQl(n)).$$

First assume that $\bF$ is separably closed.  Then since $Y_{sing}$ is a proper
closed subset of $Y$, the
groups on the
ends are zero.  By Poincar\'e duality for $Y-Y_{sing}$, we
have:

$$H_c^{2n}(Y-Y_{sing},\bQl(n))\cong H^0(Y-
Y_{sing},\bQl)^*,$$

and this last group is easily seen to be isomorphic to the
$\bQl$-dual
of

$CH^0(Y)\otimes \bQl.$  If $\bF$ is finite, then we have an exact sequence:

$$0\to H^1(\bF,H^{2n-1}(Y_{{\overline \bF}},\bQl(n)))\to H^{2n}(Y,\bQl(n))\to
H^{2n}(Y_{{\overline \bF}},\bQl(n))^{Gal({\overline\bF}/\bF)}\to 0.$$

By Deligne's theorem ([D2], Theorem 3.3.1), the weights of the geometric
Frobenius acting on $H^{2n-1}(Y_{{\overline\bF}},\bQl(n))$ are $\leq -1$ and
hence are nonzero.  This implies that
the group on the left is zero.  Using what we have just proved for separably
closed
residue fields, we get the identification ($\dagger$).  We leave the rest to
the reader.\\

{\bf Remarks 3.5}:  (i)  If the field $k$ is perfect, then
there is no
need to pass to a finite extension in Proposition 3.2.\\

(ii) It is expected that, at least after tensoring with
$\bQ$, a cycle of arbitrary codimension which is homologically equivalent to
zero on the
generic fibre
$X$ can be lifted to one which is homologically equivalent
to zero on
any regular proper model $\sX$.  \\

{\bf Proposition 3.6}:  Let $k$ be a function field in one
variable over a finite field, $X$ a smooth, projective
variety of dimension $n$ over $k$ and $l$ a prime number
different from char. $k$.  Put
$$M=H^{2n-2}(\overline X,\bQ(n)).$$  Let $S$ be a finite set
of places
of $k$ including the bad reduction places of $X$.  Then if
$z$ is a zero-cycle in the Albanese kernel $T(X)\otimes
\bQ$,

$d_n^2(z)$ is contained in the image of $$Ker[H^2(G_S,M)\to \bigoplus_{v\in
S}^{}
H^2(k_v,M)]$$

in $H^2(k,M)$,
where $d_n^2$ is the map defined in Section 2.\\

{\bf Proof}:
Let $U=Spec \: {\cal O}_{k,S}$.  Then $X$ extends to a
smooth, proper
scheme ${\cal X}_U$ over $U$ and we can define a cycle map:
$$CH^n({\cal X}_U)\to H^{2n}({\cal X}_U,\bQl(n)).$$
Consider the Leray spectral sequence for the morphism
$f:\,\sX_U\to U$:
$$H^p(U,R^qf_*\bQl(n)))\Longrightarrow H^{p+q}({\cal
X}_U,\bQl(n)).$$
Here we denote by $R^qf_*\bQl(n))$ the $l$-adic sheaf
obtained from taking the
higher direct images with finite coefficients.
Given $z$ we can extend it to a cycle on ${\cal X}_U$ of
codimension $n$ which is of
degree zero on each fibre.  This gives a cohomology class
$\cal Z$ in
$$F_1=Ker[H^{2n}({\cal X}_U,\bQl(n))\to
H^0(U,R^{2n}f_*\bQl(n))].$$  From the
spectral sequence there is a map:
$$F_1\to H^1(U,R^{2n-1}f_*\bQl(n)).$$
Claim:  $\cal Z$ is zero in $H^1(U,R^{2n-1}f_*\bQl(n)))$.
To prove the claim, we note
that since $z$ is in the Albanese kernel, its image in
$H^1(k,H^{2n-1}(\overline X,\bQl(n)))$ is zero.  Then our
claim follows from the following lemma:\\

{\bf Lemma 3.7}:  The map:
$$H^1(U,R^{2n-1}f_*\bQl(n)))\to H^1(k,H^{2n-1}(\overline
X,\bQl(n)))$$
is injective.\\

{\bf Proof}:  By the smooth and proper base change theorem,
the $l$-adic
sheaf ${\cal F}=R^{2n-1}f_*\bQl(n)$ is locally constant on
$U$.  Hence
$H^1_v(U,{\cal F})$=0 for any closed point $v$ of $U$ and
this gives the injectivity of the map.\\

	From the spectral sequence and the vanishing of the
cohomology class
$\cal Z$ in $H^1(U,R^{2n-1}f_*\bQl(n))$, we see that this
class comes from
$H^2(U,R^{2n-2}f_*\bQl(n))=H^2(G_S,H^{2n-2}(\overline
X,\bQl(n)))$.  Our task
is to show that this cohomology class vanishes in
$H^2(k_v,H^{2n-2}(\overline X,\bQl(n)))$ for all $v$ in $S$.\\

Let $\sX_v$ be a regular proper model of $X_v$ over ${\cal
O}_v$.  By Proposition 3.2, after a finite extension of $L/k_v$, we
can extend our original cycle $z$ on $X_v$ to an element of $CH^n({\cal
X}_v)\otimes \bQl$
in such a way that
its cohomology class lies in
$$H^{2n}(\sX_v,\bQl(n))^0=Ker[H^{2n}(\sX_v,\bQl(n)\to
H^0(\sO_v,R^{2n}f_*\bQl(n))].$$  Since the natural map:

$$H^2(k_v,M)\to H^2(L,M)$$

is injective for any $\bQl$-vector space $M$, it suffices to prove that our
cohomology class vanishes in
$H^2(L,M)$.

{\bf Lemma 3.8}:  $$H^{2n}(\sX_v,\bQl(n))^0=0.$$

This will prove the proposition.\\

{\bf Proof}:  From the Leray spectral sequence for $f:\,{\cal X}_v\to
\mbox{Spec}\,\sO _v$, we
get the exact sequence:
$$0\to H^2({\cal O}_v,R^{2n-2}f_*\bQl(n)))\to H^{2n}({\cal
X}_v,\bQl(n))^0\to H^1({\cal O}_v,R^{2n-1}f_*\bQl(n)))\to
0.$$

We claim that all three terms of this exact sequence are
zero.
The group on the left is zero because ${\cal O}_v$ has
cohomological
dimension one for torsion modules, hence for $l$-adic
sheaves.  As for the group
on the right, two applications of the proper base change
theorem give that it
is isomorphic to $H^1(\bF,H^{2n-1}(\overline Y,\bQl(n)))$,
where $\overline Y$ is
the geometric special fibre of ${\cal X}_v$.  By Deligne's
theorem ([D2], 3.3.1), the geometric Frobenius acts on
$H^{2n-1}(\overline Y,\bQl(n))$ with weights
which are $\leq -1$, hence nonzero.  This shows that the
group on the right
vanishes and proves the lemma.\\

\section{A Hasse Principle and Completion of the Proof}

        In this section we prove a ``Hasse'' principle in
the function
field case.  In the case where Frobenius acts semi-simply on
$H^{2n-2}(\Xb,\bQl(n))$, this is proven in Jannsen's
paper ([J2], Theorem 4 and Remark 7).
The more general result given here was shown to me by
Jannsen, and I am grateful to him for allowing me to give
his proof in this paper.  To avoid copying his earlier
proof, we
assume some familiarity with his argument.\\

{\bf Theorem 4.1} (Jannsen):  Let $k$ be a function field in
one variable over a finite field $\bF$ of characteristic $p$
and $l$ a prime number different from $p$.  Let
$V$ be a $\bQl$-linear Galois representation of
$G=Gal(k_{sep}/k)$ which is pure of weight -2, and let $S$
be a finite nonempty set of places of $k$ including those
for which the restriction to the decomposition group $G_v$
is ramified.  Then the natural map:

$$H^2(G_S,V)\to \bigoplus_{v\in S}H^2(G_v,V)$$

\noindent is injective.\\

\noindent {\bf Proof}:  Let $Y$ be a smooth, projective
model of $k$ and
let $U=Y-S, \overline U=U\times_{\bF}{\overline {\bF}}$ and
$j:\,U\to Y$
the inclusion.  Then by a theorem of Deligne ([D2], 3.4.1(iii)), $V$ is a
semi-simple $\pi_1(\overline U)$-representation.  Then we can write
$V=V_1\bigoplus V_2$,
where $\pi_1(\overline U)$
acts trivially on $V_2$ and $V_1^{\pi_1(\overline
U)}=V_{1\:{\pi_1(\overline U)}}=0$. By Jannsen's proof in
[J2], the result holds for $V_1$, so we prove it for $V_2$.
Let $F$ be the smooth $\bQl$-sheaf corresponding to $V_2$.
Then we have the exact sequence of cohomology groups:

$$(*)\,...H^2_S(Y,j_*F)\stackrel{g}{\to}H^2(Y,j_*F)\to
H^2(U,F)\stackrel{f}{\to} H^3_S(Y,j_*F)\to ...$$

By excision and the fact that
$H^2(\sO_v,j_*F)=H^2(\bF,(j_*F)_{\bF})=0$, we have:

$$H^2(k_v,F)=H^3_v(\sO_v,j_*F)$$

for all $v\in S$,
and the map $f$ is the localization map.  Thus an element of
the kernel
of the localization map comes from $H^2(Y,j_*F).$  It will
then suffice
to show that the map $g$ in the exact sequence (*) above is
surjective.

  Let $\Gamma =Gal({\overline \bF}/\bF)$ and consider the
exact sequence of $\Gamma$-modules:

$$0\to H^1(\Gamma,H^1({\overline Y},j_*F))\to H^2(Y,j_*F)\to
H^2({\overline Y},j_*F)^{\Gamma}\to 0.$$

By the fundamental theorem of Deligne ([D2], 3.2.3), the
$\Gamma$-module $H^1({\overline Y},j_*F)$ is pure of weight
-1, hence the left group is zero.  Consider the following
diagram:

$$\matrix{H^2_S(Y,j_*F)&\to&H^2_{\overline S}({\overline
Y},j_*F)^{\Gamma}\cr
\downarrow&&\downarrow\cr
H^2(Y,j_*F)&\to&H^2({\overline Y},j_*F))^{\Gamma}}.$$

We have just shown that the bottom horizontal map is an
isomorphism and our goal
is to show that the left vertical arrow is surjective.
Since the top horizontal map is surjective, it will be
enough to prove the surjectivity of the right vertical map.
Let $\tilde{\sO_v}$ denote the strict henselization of
$\sO_v$ and $\tilde{k_v}$ its fraction field.  Then we have:

$$ H^1(\tilde{k_v},F)\cong H^2_v(\tilde{\sO_v},j_*F).$$

Now the pro-$l$-part of the absolute Galois group of
$\tilde{k}$ is isomorphic to $\bZl(1)$ as a $\Gamma$-module.
Then by Pontryagin duality and the fact that
$\pi_1({\overline U})$ is
assumed to act trivially on $V_2$, we have

$$H^1(\tilde{k_v},F)^{\Gamma}=(\hat{V_2}(1)_{\Gamma})^*,$$
where * denotes $\bQl$-linear dual and $\hat{V_2}$ denotes
the dual representation.  But by Poincar\'e duality, we
have:

$$H^2({\overline Y},j_*F)^{\Gamma}\cong
(\hat{V_2}(1)_{\Gamma})^*,$$

and since $S$ is nonempty, we obviously have injectivity of
the dual map:

$$(\hat{V_2}(1)_{\Gamma})\to [\bigoplus_{v\in
S}\hat{V_2}(1)_{\Gamma}]^.$$  This completes the proof of
the Theorem.

	Now we can prove the main theorem of this paper:\\

{\bf Theorem 4.2}:  Let $X$ be a smooth, projective,
geometrically connected
variety of dimension $n$ over a function field $k$ in one
variable over a finite field of characteristic $p$, and
$l$ a prime number different from $p$.  Assume Hypothesis
3.1 above.
Then the map $d_n^2$ of
section 1 is zero. \\

{\bf Proof}:  Combine Proposition 3.6 and Theorem 4.1, taking
$V=H^{2n-2}(\Xb,\bQl(n))$.\\

{\bf Theorem 4.3}:  Let $X$ be a smooth, projective,
geometrically connected
surface over a function field in one variable over a finite
field of characteristic
$p\geq 5$.  Then the map $d_2^2$ is zero.\\

{\bf Proof}:  Abhyankar [A] has proved resolution of
singularities for 3-folds
over an algebraically closed field of characteristic greater than or equal
to 5.
Given $X$, we can pass to a finite extension over which it
has a smooth
proper model.  Then the kernel of $d_2^2$ will be killed by
the degree
of this extension, hence it is zero.

\newpage
\section*{Appendix}

	In this appendix we prove a technical lemma which shows that the
Abel-Jacobi map as defined via the Hochschild-Serre spectral sequence
agrees with the map defined by Kummer theory on $Alb(X)$.  This
fact is probably known to several people, but I have been unable to find
a reference.\\

{\bf Lemma}:  Let $X$ be a smooth projective geometrically connected
variety of dimension $n$ over a field $k$, $l$ a prime number different from
$\mbox{char.}k$.  Consider the map $d_n^1$ defined in \S 1 above and the
map:
$$f_n^1:  A_0(X)\to H^1(k,H^{2n-1}(\Xb,\bQl(n)))$$

defined as the composite of the Albanese map

$$\alpha:  A_0(X)\to Alb(X)$$

and the map
$$Alb(X)\to H^1(k,H^{2n-1}(\Xb,\bQl(n)))$$

obtained from Kummer theory on $Alb(X)$.  Then $d_n^1=f_n^1$.\\

{\bf Proof}:  First assume that $X$ is a curve.  The group $A_0(X)$ is
generated by elements of the form
$\mbox{cor}_{L/k}(P-Q)$, where $L$ is a finite extension of $k, P$ and $Q$ are
rational points of $X_L$ and
$$\mbox{cor}:\,A_0(X_L)\to A_0(X)$$
denotes corestriction.  Then it suffices to prove that for any such $L/k$ and
any two
$L$-rational points $P,Q$ of $X$, the two maps above give the same
element in $H^1(L,H^{1}(\Xb,\bZl(1)))$.  We suppress $L$ in what
follows.  Let $z=P-Q$, let $Z$ be its support and $U=X-Z$.
Taking cohomology with support in ${\overline Z}$, we
easily get the exact sequence:

$$(*):\, 0\to H^1(\Xb,\bZ/l^m(1))\to H^1({\overline U},\bZ/l^m(1))\to
\bZ/l^m\to 0,$$

where the last group is generated by the class of $z$.  Jannsen
has shown ([J5], Lemma 9.4) that the extension given by (*) gives the element
of
$H^1(k,H^1(\Xb,\bZ/l^m(1))$ which is associated to $P-Q$ via the map
$d_1^1$.  Now
consider the Kummer sequence for $Pic^0(\Xb)$:

$$0\to H^1(\Xb,\bZ/l^m(1))\to Pic^0(\Xb)\stackrel{l^m}{\to}Pic^0(\Xb)\to
0,$$
where the group on the left is identified with the $l^m$-torsion of
$Pic^0(\Xb)$.
As is well-known, the boundary map:
$$Pic^0(\Xb)^G\to H^1(k,H^1(\Xb,\bZ/l^m(1)))$$
in the $G$-cohomology of this
sequence may be may be obtained by pulling back the Kummer sequence by the
$G$-map
$$\bZ\to Pic^0(\Xb)$$

corresponding to an element $y$ of $Pic^0(\Xb)^G$.  One easily sees that the
pullback extension is:

$$(**):\: 0\to H^1(\Xb,\bZ/l^m(1))\to D\to \bZ\to 0,$$

where $D=\{x\in Pic^0(\Xb): l^mx=ny\, \mbox{for some}\, n\in \bZ\}$.  We may
define a $G$-map

$$D/l^m\to H^1({\overline U},\bZ/l^m(1))$$

as follows:  Let $x\in D$, let ${\cal L}$ be the isomorphism class of line
bundle
associated to $x$ and ${\cal M}$ the isomorphism class of
line bundle associated to $y$.   Choose an isomorphism
$${\cal L}^{\otimes l^m}\to {\cal M}^{\otimes n}.$$  Such an isomorphism is
unique modulo $l^m$, since if we had two such, then we would get an isomorphism
from
${\cal O}_{\Xb}$ to itself.  Since $\Xb$ is proper, this is given by an
element of $\kb^*$, and hence must be the
identity modulo $l^m$.  But this means that the two isomorphisms were the
same modulo $l^m$.  Since the class of $y$ is
supported on ${\overline Z}$, the restriction of this isomorphism to
${\overline U}$ gives a
trivialization:

$${\cal L}^{\otimes l^m}\cong {\cal O}_{\overline U},$$
and this data defines in a canonical way an element of $H^1({\overline
U},\bZ/l^m(1))$.  It is easy to verify that via this map, the two
extensions (*) and (**) give the same element of $H^1(k,H^1(\Xb,\bZ/l^m(1)))$.
Note that for this argument, we have used $\bZl$-coefficients.

	In the case of $X$ of arbitrary dimension, Poincar\'e duality and the Weil
pairing give an identification:

$$H^{2n-1}(\Xb,\bQl(n))=V_l(Alb(\Xb)).$$

Now given a zero-cycle
$z$, if $k$ is perfect, we can find a smooth, projective curve $C$ going
through
it ([AK], [C]).  By
functoriality properties of the maps in question, we can easily reduce
to the case of $C$ which we have done above.  If $k$ is not perfect,
such a curve $C$ may be found after an inseparable extension of $k$.
Then an easy restriction-corestriction argument shows that the maps are
the same.  This completes the proof
of the lemma.

\newpage
\section*{References}
[A]  Abhyankar, S.:  {\bf Resolution of Singularities of
Embedded Algebraic Surfaces}, Academic Press 1966\\

[AK]  Altman, A. and Kleiman, S.  {\em Bertini theorems for
hypersurface
sections containing a subscheme}, Comm. in Algebra {\bf 7}
(1979),
775-790\\

[Be1]  Beauville, A. {\em Quelques remarques sur la transformation de Fourier
dans l'anneau de Chow d'une vari\'et\'e ab\'elienne}, in {\bf Algebraic
Geometry:  Tokyo/Kyoto 1982}, Lecture Notes in Mathematics 1016, Springer
1983\\

[Be2]  Beauville, A. {\em Sur l'anneau de Chow d'une vari\'et\'e
ab\'elienne}, Math. Ann. {\bf 273} (1986), 647-651\\

[B1]  Beilinson, A. {\em Height pairing on algebraic
cycles}, in {\bf Current Trends in Arithmetical Algebraic
Geometry}, K. Ribet editor, Contemporary Mathematics Volume
67, American Mathematical Society 1987\\

[B2]  Beilinson, A. {\em Absolute Hodge cohomology}, in {\bf
Applications of Algebraic K-Theory to Algebraic Geometry and
Number
Theory}, S. Bloch, R.K. Dennis, E.M. Friedlander and M.
Stein editors,
Contemporary Mathematics Volume 55, American Mathematical
Society 1986\\

[Bl1]  Bloch, S. {\em Algebraic cycles and values of $L$-
functions}, J.
f\"ur die Reine und Ang. Math. {\bf 350} (1984), 94-108\\

[Bl2]  Bloch, S. {\em $K_2$ of Artinian $\bQ$-algebras with
application
to algebraic cycles}, Communications in Algebra {\bf 3}
(1975), 405-428\\

[Bl3]  Bloch, S. {\em Some elementary theorems about algebraic cycles on
abelian varieties}, Inventiones Math. {\bf 37} (1976), 215-228\\

[Bl4]  Bloch, S.  {\em Lectures on Algebraic Cycles}, Duke Univ. Press, 1980\\

[BK]  Bloch, S. and Kato, K.  {\em L-functions and Tamagawa
numbers of
motives}, in {\bf Grothendieck Festschrift}, Volume I, P.
Cartier
editor, Progress in Mathematics, Birkh\"auser 1990\\

[BLR]  Bosch, S., L\"ukebohmert, W. and Raynaud, M.  {\bf N\'eron Models},
Springer 1990\\

[CT-O]  Colliot-Th\'el\`ene, J.-L. and Ojanguren, M. {\em Vari\'et\'es
unirationnelles nonrationnelles:  au d\'ela de l'exemple d'Artin-Mumford},
Inventiones Math. {\bf 97} (1989), 141-158\\

[CT-R]  Colliot-Th\'el\`ene, J.-L. and Raskind, W.  {\em
Groupe de Chow
de codimension deux des vari\'et\'es d\'efinies sur un corps
de nombres:
un th\'eor\`eme de finitude pour la torsion}, Inventiones
Math. {\bf
105} (1991), 221-245\\

[C]   Coray, D. {\em Two remarks on the Bertini theorems},
unpublished
typescript, 1980\\

[D1]  Deligne, P. {\em  Th\'eor\`eme de Lefschetz difficile
et crit\`eres de d\'eg\'en\'eresence de suites spectrales}
Publ. Math. I.H.E.S. {\bf 35} (1968), 107-126\\

[D2]  Deligne, P.  {\em La Conjecture de Weil II}, Publ.
Math. I.H.E.S. {\bf 52}
(1980), 137-252\\

[E]  Ekedahl, T.  {\em On the adic formalism}, in {\bf
Grothendieck
Festschrift}, Volume II, P.
Cartier editor, Progress in Mathematics, Birkh\"auser 1990\\

[Gi]  Gillet, H.  {\em Gersten's conjecture for the K-theory
with torsion
coefficients of a discrete valuation ring}, J. of Algebra
{\bf 103}
(1986), 377-380\\

[J1]  Jannsen, U. {\em Continuous \'etale cohomology}, Math.
Annalen {\bf 280}
(1988), 207-245\\

[J2]  Jannsen, U. {\em On the $l$-adic cohomology of
varieties over number
fields and its Galois cohomology}, in {\bf Galois Groups
over Q}, Y. Ihara,
K. Ribet and J.-P. Serre editors, Springer 1989\\

[J3]  Jannsen, U. {\em letter to B. Gross}, 1989\\

[J4]  Jannsen, U. {\em Motivic sheaves and filtrations on
Chow groups}, in {\bf Motives:  Seattle 1991},  U. Jannsen,
S. Kleiman  editors\\

[J5]  Jannsen, U. {\bf Mixed Motives and Algebraic K-Theory}, Lecture
Notes in Mathematics 1400, Springer 1990\\

[L]   Lang, S. {\bf Diophantine Geometry}, Springer 1983\\

[M1]  Milne J.S.  {\bf Etale Cohomology}, Princeton
Mathematical Series Volume 33, Princeton University
Press 1980\\

[M2]  Milne, J.S. {\em Zero-cycles on algebraic varieties in
non-zero
characteristic:  Rojtman's theorem}, Composition Math. {\bf
47} (1982), 271-287\\

[Mu]  Murre, J. {\em On a conjectural filtration on the Chow
groups of an algebraic variety}, preprint 1992\\

[N]  Nori, M.  {\em Algebraic cycles and Hodge-theoretic connectivity},
Inventiones Math. {\bf 105} (1993) 349-373\\

[Ro]  Roitman, A.A.  {\em The torsion of the group of 0-cycles modulo
rational equivalence},  Annals
of Math. {\bf 111} (1980), 553-569\\

[S1]  Saito, S. {\em Motives and filtrations on Chow
groups}, preprint 1993\\

[S2]  Saito, S.  {\em On the cycle map for torsion algebraic
cycles of
codimension two}, Inventiones Math. {\bf 106} (1991), 443-460\\

[S3]  Saito, S. {\em Remarks on algebraic cycles}, preprint 1993\\

[SH] Saito, H. {\em Generalization of Abel's Theorem and
some finiteness
properties of zero-cycles on surfaces}, preprint\\

[T]  Thomason, R.  {\em Absolute Cohomological Purity},
Bull. Math. Soc. France {\bf 112} (1984), 397-406\\

\end{document}